\documentclass[9pt]{article}
\usepackage[russian]{babel}
\usepackage{amssymb}
\usepackage{amsthm}
\usepackage{bbm}
\usepackage{amsmath}

\newtheorem{lemma}{Лемма}
\newtheorem{_theorem1}{Theorem}
\newtheorem{_theorem}{Теорема}
\newtheorem{_cor1}{Corollary}

\begin{document}

\section*{Bykovskii's theorem and generalized Larcher's theorem}
\centerline{Dmitry M. Ushanov\footnote{research is supported by RFBR grant 12-01-00681}}
\subsection*{Summary}

For $a = (a_1, \dots, a_s) \in \mathbb{Z}^s$ denote by
$$
\mathcal{K}(a) =
\left\{
	\left(
		\left\{
		\frac{a_1 k}{N}
		\right\},
		\cdots,
		\left\{
		\frac{a_s k}{N}
		\right\}
	\right)
	\mid
	k = 1,2,\dots,N
\right\}
$$
the sequence of \textit{Korobov lattice points} (see \cite{Kor}).

For any sequence 
$\Xi = \left\{ \xi_x \in [0,1]^s,\,\,\,\, x=0,1,2,...,N-1 \right\}$
denote by
$$
D(\Xi) = \sup_{\gamma_1,\dots\,\gamma_s \in [0,1]} |\#\{x : 0\leq x < N,\ \xi_x\in [0,\gamma_1)\times\cdots\times[0,\gamma_s)\} - \gamma_1\cdots\gamma_s N|
$$
the \textit{discrepancy} of the sequence.

G. Larcher proved the following theorem (see \cite{Lar}).

\textbf{Theorem A.} \textit{
There is a constant $c,$
such that for every integer
$N$ there exists an integer $g,$ with $(g, N) = 1$ and for
the sequence 
$\mathcal{K}(1,g)$ is
$$
D_N(\mathcal{K}(1,g)) < c \frac{N\log{N} \log\log{N}}{\phi(N)}.
$$}

Recently V. Bykovskii proved the following result
(see \cite{By2}, \cite{By3}).

\textbf{Theorem B.} \textit{
For every integer
$s \geqslant 2$
there is a constant $c,$
such that for every integer $N$
there exists an $a = (a_1,\dots,a_s) \in \mathbb{Z}^s$ with
$$
D_N(\mathcal{K}(a_1,\dots,a_s)) < c \log^{s-1}{N} \log\log{N}.
$$}

Also, G. Larcher obtained the following result
(see \cite{Lar}, Corollary 5).
For every $k\in \mathbb{N}$ there is a constant $c(k)$ such that for every $N\in \mathbb{N},$ there exists an $x\in \mathbb{N}$ with $(x,N)=1,$
and
$$
\sum {b_i} < c(k) ( \log N \log\log N )^2,
$$
where $x^k/N = [b; b_1, b_2, \dots, b_l]$ is a continued
fraction expansion.
Therefore (see \cite{Ostr}), one has
$$
D(\mathcal{K}(1,x^k)) \ll ( \log N \log\log N )^2.
$$

N. Moshchevitin and D. Ushanov \cite{UshMo} improved this result
of G. Larcher and proved the following theorem.

\textbf{Theorem C.} \textit{ Let $p$ be prime, $U$ be a multiplicative  subgroup in $\mathbbm{Z}_p^*.$  For  $v \neq 0$ we consider the set $R = v\cdot U$ and let
$$
\#R \geq 10^8 p^{7/8}\log^{5/2}{p}.
$$
Then there exists an element $a \in R,$ $a/p=[b_1, b_2, \cdots, b_l],$ $b_i = b_i(a),$ $l = l(a)$ with
$$
\sum_{i=1}^{l} b_i \leq 500\log{p}\log\log{p}.
$$
}

We would like to note that recently Mei-Chu Chang \cite{Chang} independently obtained Theorem C.

In present paper we use results of Bykovskii and obtain the following
generalization of Larher's theorem.

\begin{_theorem1}
Let
$$
\alpha _m = \frac{2^{m-1}}{2^{m+1}-m-2}, \ \ \ \beta _m = \frac{2^m}{2^{m+2}-m-4},
$$
$$
n_\delta = \max( 3, \left\lceil \frac{81}{\delta} + \frac{1}{\ln 2} - 160 \right\rceil).
$$
Define $s'(\delta)$ and $s''(\delta)$ as follows
$$
s'(\delta) =
\begin{cases}
\left\lceil \frac{2 m^2 \left(\delta{} - 1\right)}{4\delta(2^{-m} - 1) + 1}\right\rceil + 2,
	& \text{if $\beta_m \leqslant \delta < \alpha_m$} \\
\left\lceil \frac{2 m \left(\delta - 1\right)\left(m + 1\right)}{\delta ( 3 \cdot 2^{-m} - 4 ) + 1}\right\rceil + 2,
	& \text{if $\alpha_{m+1} \leqslant \delta < \beta_m$}, \\
\end{cases}
$$
$$
s''(\delta) = \left\lfloor \frac{(1-\delta)2^{n_\delta}}{(\delta \frac{160+n_\delta}{81} - 1)\cdot 0.45} \right\rfloor + 3.
$$
Set $s_{\min}(\delta) = \min(s'(\delta), s''(\delta)).$

Let $\delta \in (0,\ 1)$ and integer $s \geqslant s_{\min}(\delta).$
Let $p \geqslant 3$ be prime, $G \subset \mathbb{Z}_p^*$ be a subgroup. If
$\#G \geqslant p^\delta$
then there exist elements $a_1,\dots,a_s \in G$ with
$$D_p(\mathcal{K}(a_1,\dots,a_s)) \ll_{\delta,\ s} \log^{s-1}p\log\log p.$$
\end{_theorem1}

\begin{_cor1}
For small values of $m$ in Theorem 1 we get the following.
Let $s$ be an integer and $\delta$ be a real such that
$$
s \geqslant 
\begin{cases}
3, & \delta\in [3/4,\ 1) \\
4, & \delta\in [2/3,\ 3/4) \\
5, & \delta\in [14/23,\ 2/3) \\
6, & \delta\in [4/7,\ 14/23) \\
7, & \delta\in [6/11,\ 4/7) \\
8, & \delta\in [10/19,\ 6/11) \\
9, & \delta\in [22/43,\ 10/19) \\
10, & \delta\in [1/2,\ 22/43) \\

11, & \delta\in [17/35, 1/2) \\
\end{cases}
s \geqslant
\begin{cases}
12, & \delta\in [9/19, 17/35) \\
13, & \delta\in [19/41, 9/19) \\
14, & \delta\in [5/11, 19/41) \\
15, & \delta\in [21/47, 5/11) \\
16, & \delta\in [11/25, 21/47) \\
17, & \delta\in [23/53, 11/25) \\
18, & \delta\in [3/7, 23/53) \\
19, & \delta\in [25/59, 3/7) \\
20, & \delta\in [13/31, 25/59). \\
\end{cases} 
$$
Let $p \geqslant 3$ be prime, $G \subset \mathbb{Z}_p^*$ be a subgroup. If
$\#G \geqslant p^\delta$
then there exist elements $a_1,\dots,a_s \in G$ with
$$D_p(\mathcal{K}(a_1,\dots,a_s)) \ll_{\delta,\ s} \log^{s-1}p\log\log p.$$ 
\end{_cor1}

The author is grateful to Prof. Igor Shparlinski for pointing out an opportunity of improvement of the main result of the paper.

\newpage

\subsection*{1. Введение}

Пусть $a = (a_1, \dots, a_s)$ --- произвольный набор целых чисел,
а $N \geqslant 3$ --- натуральное число. \textit{Сеткой Коробова}
(см. \cite{Kor})
называется множество
$$
\mathcal{K}(a) =
\left\{
	\left(
		\left\{
		\frac{a_1 k}{N}
		\right\},
		\cdots,
		\left\{
		\frac{a_s k}{N}
		\right\}
	\right)
	\mid
	k = 1,2,\dots,N
\right\}.
$$
Рассмотрим решетку
$$
\Gamma_N(a) =
\left\{
	(m_1,\dots,m_s) \in \mathbb{Z}^s \mid a_1 m_1 +\cdots+ a_s m_s \equiv 0 \pmod{N}
\right\},
$$
несложно проверить, что 
$$
\mathcal{K}(a) = \Gamma_N^*(a) \cap [0,1)^s,
$$
где $\Gamma_N^*(a)$ --- двойственная к $\Gamma_N(a)$ решетка.

Для последовательности
$\Xi = \left\{ \xi_x \in [0,1]^s,\,\,\,\, x=0,1,2,...,N-1 \right\}$
\textit{отклонением} называется следующая величина:
$$
D(\Xi) = \sup_{\gamma_1,\dots\,\gamma_s \in [0,1]} |\#\{x : 0\leq x < N,\ \xi_x\in [0,\gamma_1)\times\cdots\times[0,\gamma_s)\} - \gamma_1\cdots\gamma_s N|.
$$

В \cite{Lar} Г. Ларчер получил следующий результат.

\textbf{Теорема A.} \textit{
Существует постоянная $c,$
такая, что для любого
натурального $N$ существует натуральное $g,$ $(g, N) = 1$ и отклонение сетки
$\mathcal{K}(1,g)$ оценивается
$$
D_N(\mathcal{K}(1,g)) < c \frac{N\log{N} \log\log{N}}{\phi(N)}.
$$}

Доказательство теоремы элементарно и опирается на аппарат цепных дробей.

Недавно В.А. Быковский получил следующий результат (см. \cite{By2}, \cite{By3}).

\textbf{Теорема B.} \textit{
Для любого натурального $s \geqslant 2$
существует постоянная $c,$
такая, что для любого натурального $N$ найдется такой набор
$(a_1,\dots,a_s),$ что отклонение сетки $\mathcal{K}(a_1,\dots,a_s)$ оценивается
$$
D_N(\mathcal{K}(a_1,\dots,a_s)) < c \log^{s-1}{N} \log\log{N}.
$$}

Доказательство этой теоремы использует аналитический аппарат.

Ларчер получил также следующий результат
для вычетов степени $k$ (см. \cite{Lar}, Следствие 5).
Для любого $k\in \mathbb{N}$ существует константа $c(k)$ такая,
что для любого $N\in \mathbb{N}$ найдется $x\in \mathbb{N},$ $(x,N)=1,$
$1 \leqslant x < N,$ такой что существует $x_1,$ $x \equiv x_1^k \pmod{N}$
и для разложения в цепную дробь $\frac{x}{N} = [0;b_1,\dots,b_l]$
выполнено
$$
\sum {b_i} < c(k) ( \log N \log\log N )^2.
$$
Отсюда следует (см. Островский \cite{Ostr}) оценка на отклонение двумерной
последовательности: для такого $x$ будет верно
$$
D(\mathcal{K}(1,x)) \ll ( \log N \log\log N )^2.
$$

В работе \cite{UshMo} получено усиление теоремы Ларчера
на случай когда вместо множества вычетов фиксированной степени
при простом $N$ рассматривается
произвольная достаточно большая мультипликативная
подгруппа группы $\mathbb{Z}_p^*.$

\textbf{Теорема C.} \textit{ Пусть $p$ простое, $G$ 
--- мультипликативня подгруппа в $\mathbbm{Z}_p^*.$  Пусть
$$
\#G \geq 10^8 p^{7/8}\log^{5/2}{p}.
$$
Тогда существует элемент $a \in G,$ $a/p=[0; b_1, b_2, \cdots, b_l],$ с
$$
D(\mathcal{K}(1,a)) \ll \sum_{i=1}^{l} b_i \leq 500\log{p}\log\log{p}.
$$}

Пользуясь результатами В.А. Быковского (см. \cite{By2}, \cite{By3}),
в настоящей работе
мы обобщили результат Ларчера на большие размерности
и получили следующую теорему.

\begin{_theorem}
\label{sdf2234}
Обозначим
$$
\alpha _m = \frac{2^{m-1}}{2^{m+1}-m-2}, \ \ \ \beta _m = \frac{2^m}{2^{m+2}-m-4},\ \ \ 
n_\delta = \max\left( 3, \left\lceil \frac{81}{\delta} + \frac{1}{\ln 2} - 160 \right\rceil\right).
$$
Определеним функции $s'(\delta)$ и $s''(\delta)$ следующим образом
$$
s'(\delta) =
\begin{cases}
\left\lceil \frac{2 m^2 \left(\delta{} - 1\right)}{4\delta(2^{-m} - 1) + 1}\right\rceil + 2,
	& \text{если $\beta_m \leqslant \delta < \alpha_m$}, \\
\left\lceil \frac{2 m \left(\delta - 1\right)\left(m + 1\right)}{\delta ( 3 \cdot 2^{-m} - 4 ) + 1}\right\rceil + 2,
	& \text{если $\alpha_{m+1} \leqslant \delta < \beta_m$}, \\
\end{cases}
$$
$$
s''(\delta) = \left\lfloor \frac{(1-\delta)2^{n_\delta}}{(\delta \frac{160+n_\delta}{81} - 1)\cdot 0.45} \right\rfloor + 3.
$$
Пусть $s_{\min}(\delta) = \min(s'(\delta), s''(\delta)).$

Если $s \geqslant 3$ и $\delta \in (0,\ 1)$ таковы, что $s \geqslant s_{\min}(\delta),$ то
для любого простого $p \geqslant 3,$ любой подгруппы $G$ мультипликативной
группы $\mathbb{Z}_p^*,$ такой, что
$\#G \geqslant p^\delta,$
найдутся элементы $a_1,\dots,a_s \in G,$ для которых верно
$$D_p(\mathcal{K}(a_1,\dots,a_s)) \ll_{\delta,\ s} \log^{s-1}p\log\log p.$$
\end{_theorem}

\subsection*{2. Основые идеи доказательства}

Всюду далее мы будем использовать следующие обозначения.
Пусть $G$ --- подгруппа мультипликативной группы $\mathbb{Z}_p^*$ чисел по модулю $p$.
Обозначим через $A(P_1,...,P_s)$ число решений
сравнения
$$a_1m_1 + \cdots + a_s m_s \equiv 0 \pmod{p}$$
в целых $a_1,\dots,a_s \in G$ и целых $m_1,\dots,m_s,$ таких, что
$$
1\leqslant m_1 \leqslant P_1, \dots, 1\leqslant m_s \leqslant P_s.
$$

\begin{lemma}
\label{main_lemma_3}
Пусть параметры $s$ и $\delta$ выбраны так, как в теореме \ref{sdf2234}.
Тогда
$$A(P_1,\dots,P_s) \ll \frac{\#G ^ s}{p} P_1 \cdots P_s.$$
\end{lemma}

\begin{proof}
Обозначим
$e_p(t) = \exp( 2\pi i t / p ),$
$$S(t, G) = \sum_{g\in G} e_p(t g),$$
$$S(G) = \max_{t \neq 0} \left| S(t, G) \right|.$$

Сначала рассмотрим случай, когда $s \geqslant s'(\delta).$

Заметим, что количество решений сравнения равно
$$A(P_1,\dots,P_s) = \frac{1}{p}
  \sum_{1 \leqslant m_i \leqslant P_i}
  \sum_{a_i \in G}
  \sum_{n=1}^{p}
   e_p(n(a_1 m_1 + \dots + a_s m_s)).
$$
Воспользовавшись неравенством Коши-Буняковского, получим
\begin{equation}
\label{AEst}
A(P_1,\dots,P_s) \leqslant
 \frac{1}{p} \sum_m S^{s-2}(G) \cdot p \#G.
\end{equation}

Для оценки $S(G)$ сверху мы воспользуемся результатом Конягина (см. \cite{Kng}). А именно, в обозначениях формулировки
теоремы \ref{sdf2234}
\begin{enumerate}
\item если существует $m\in \mathbb{N},$ что $\beta_m \leqslant \delta < \alpha_m$, то
	$$S(G) \ll_m \#G^{1 - \frac{2}{m^2} + \frac{1}{2^{m-1}m^2}} p^{\frac{1}{2m^2}},$$
\item если же существует $m\in \mathbb{N},$ что $\alpha_{m+1} \leqslant \delta < \beta_m,$ то
	$$S(G) \ll_m \#G^{1 - \frac{2}{m(m+1)} + \frac{3}{2^{m+1}m(m+1)}} p^\frac{1}{2m(m+1)}.$$
\end{enumerate}

Подставляя оценки на $S(G)$ в оценку (\ref{AEst}) для $A(P_1,\dots,P_s),$
мы получим
$$
A(P_1,\dots,P_s) \leqslant \frac{\#G^s}{p} P_1\cdots P_s.
$$

В случае $s \geqslant s''(\delta)$
мы рассуждаем так же,
разница заключается в другой оценке $S(G).$

Мы используем следующий результат Гараева (см. \cite{Gar}, теорема 4.1).
Если $3 \leqslant n\leqslant 1.44 \log \log p$ --- натуральное число,
$c>0$ --- произвольная постоянная, $X_1, \dots, X_n$ --- подмножества
$Z_p^*,$ удовлетворяющие условию
$$\#X_1 \cdot \#X_2 \cdot \left( \#X_3 \cdots \#X_n \right)^{1/81}> p^{1+c},$$
то имеет место оценка
$$\left| \sum_{x_1 \in X_1} \dots \sum_{x_n \in X_n} e_p(x_1\cdots x_n) \right|
\ll \#X_1 \cdots \#X_n p^\frac{-0.45c}{2^n}.$$

Взяв $X_1 = tG,$ $X_2 = G,\dots,X_n = G,$ и подобрав подходящим образом
параметры $c$ и $n,$ мы получим лемму.
\end{proof}

\subsection*{3. Схема доказательства теоремы \ref{sdf2234}}
Выберем параметр $Q = \frac{ p }{ 2\log^{s+1}p}.$
Пусть $a = (a_1, \dots, a_s).$
Для целочисленного вектора $m = (m_1, \dots, m_s)$ определим
$$
H(m) = \max(1, |m_1|)\cdots \max(1, |m_s|).
$$
Определим функцию 
$$
q(a) = \min \{ H(m) \mid a_1 m_1 + \cdots + a_s m_s \equiv 0 \pmod{p} \}.
$$
Определим множества
$\Omega_1$ и $\Omega$ следующим образом:
$$
\Omega_1 = \{a \in G^s \mid q(a) \leqslant Q\},
$$
$$
\Omega = \{a \in G^s \mid q(a) > Q\}.
$$

Очевидно, что $\#(\Omega_1 \cup \Omega) = \#G^s.$
С помощью леммы \ref{main_lemma_3} мы показываем, что $\#\Omega_1 \leqslant \#G^s /2,$ следовательно
$\#\Omega > \#G^s /2.$

Ненулевой узел $\gamma$ решетки $\Gamma$ называется
\textit{относительным минимумом}, если не найдется другого ненулевого узла $\gamma'$ из $\Gamma,$
для которого
$$
|\gamma'_1| \leqslant |\gamma_1|,\; \dots, \; |\gamma'_s| \leqslant |\gamma_s|.
$$
Обозначим через $\mathfrak{M}(\Gamma)$ множество всех относительных минимумов решетки $\Gamma.$

Далее мы воспользуемся следующим
результатом В.А. Быковского (см. \cite{By3}): для любой сетки
Коробова верно неравенство
$$
D_N(\mathcal{K}(a)) \ll N \sum_{m \in \mathfrak{M}(\Gamma_N(a))} \frac{1}{H(m)}.
$$

Таким образом, усредняя по множеству $\#\Omega$, мы получим
$$
\min_{a\in G^s} D(\mathcal{K}(a))
\ll \frac{p}{\#\Omega}
\sum_{a\in G^s}
\sum_{Q\leqslant H(m) \leqslant p}\frac{\delta_p(a_1 m_1 + \cdots + a_s m_s)}{H(m)}.
$$
Оценивая последнюю сумму, мы получаем теорему.


\vspace{40pt}
Dmitry M. Ushanov

\textit{E-mail}: ushanov.dmitry@gmail.com

\end{document}